\documentclass[opre,nonblindrev,copyedit]{informs1} % for review, not blinded
\usepackage{natbib}
 \NatBibNumeric
 \def\bibfont{\small}%
 \bibpunct[, ]{[}{]}{,}{n}{}{,}%

\usepackage[colorlinks=true,breaklinks=true,bookmarks=true,urlcolor=blue,
     citecolor=blue,linkcolor=blue,bookmarksopen=false,draft=false]{hyperref}

\def\EMAIL#1{\href{mailto:#1}{#1}}% When hyperref is used, otherwise outcomment 
\def\URL#1{\href{#1}{#1}}         % When hyperref is used, otherwise outcomment 

\TheoremsNumberedBySection  % (Theorem 1.1, Lema 1.1, Theorem 1.2)

\EquationsNumberedBySection % (1.1), (1.2), ...

\newcommand{\mb}[1]{\ensuremath{\boldsymbol{#1}}}
\usepackage{amsmath,  amssymb, epsfig, graphicx, float, url,tabularx,multirow}
\usepackage{algorithm, algorithmicx, algpseudocode,bbm,subcaption,booktabs}
\usepackage{enumitem}
\usepackage{footnote}
\makesavenoteenv{tabular}
\makesavenoteenv{table}

\newcommand\blfootnote[1]{%
  \begingroup
  \renewcommand\thefootnote{}\footnote{#1}%
  \addtocounter{footnote}{-1}%
  \endgroup
}

\begin{document}
%%%%%%%%%%%%%%%%

%\setcounter{page}{1}
%\DOI{10.1287/xxxx.0000.0000}%


\RUNAUTHOR{El Housni and Goyal}
 \RUNTITLE{Beyond Worst-case: A Probabilistic Analysis of Affine Policies in Dynamic Optimization}
%\RUNTITLE{}

% \TITLE{Bundling Information Goods of Decreasing Value}
%\TITLE{Power of Affine Policies: Optimal approximation for Budgeted Uncertainty Sets}
\TITLE{Beyond Worst-case: A Probabilistic Analysis of Affine Policies in Dynamic Optimization}
\ARTICLEAUTHORS{%
\AUTHOR{Omar El Housni}
\AFF{Dept of Industrial Engineering and Operations Research, Columbia University, New York NY 10027, \EMAIL{oe2148@columbia.edu}, \URL{http://www.columbia.edu/~oe2148/}}
\AUTHOR{Vineet Goyal}
\AFF{Dept of Industrial Engineering and Operations Research, Columbia University, New York NY 10027, \EMAIL{vg2277@columbia.edu}, \URL{http://www.columbia.edu/~vg2277/}}
} % end of the block

\ABSTRACT{%
 Affine policies (or control) are widely used as a solution approach in dynamic optimization where computing an optimal adjustable solution is usually intractable. While the worst case performance of affine policies can be significantly bad, the empirical performance is observed to be near-optimal for a large class of problem instances. For instance, in the two-stage dynamic robust optimization problem with linear covering constraints and uncertain right hand side, the worst-case approximation bound for affine policies is $O(\sqrt m)$ that is also tight (see Bertsimas and Goyal~\cite{BG10}), whereas observed empirical performance is near-optimal. In this paper, we aim to address this stark-contrast between the worst-case and the empirical performance of affine policies. In particular, we show that  with high probability affine policies give a good approximation for  two-stage dynamic robust optimization problems  on random instances generated from a large class of distributions; thereby, providing a theoretical justification of the observed empirical performance. 
The approximation bound depends on the distribution, but  it is significantly better than the worst-case bound for a large class of distributions.
 %On the other hand, we also present a distribution such that the performance bound for affine policies on instances generated according to that distribution is $\Omega(\sqrt m)$ with high probability. This demonstrates that the empirical performance of affine policies can depend on the generative model for instances.
}%

%\KEYWORDS{deterministic inventory theory; infinite linear programming duality; 
%  existence of optimal policies; semi-Markov decision process; cyclic schedule}
%\MSCCLASS{Primary: 90B05; secondary: 90C40, 90C90}
%\ORMSCLASS{Primary: Inventory/production: deterministic multi-item;
%  secondary: dynamic programming/optimal control: deterministic 
%  semi-Markov; programming: infinite dimensional}
%\HISTORY{Received November 20, 2003; revised March 8, 2004, and March 26, 2004.}

%\MSCCLASS{Primary: 90C15 , 90C47}
%\ORMSCLASS{Primary: ; secondary: }

\maketitle
%%%%%%%%%%%%%%%%%%%%%%%%%%%%%%%%%%%%%%%%%%%%%%%%%%%%%%%%%%%%%%%%%%%%%%

\section{Introduction}
\blfootnote{A preliminary version of this paper with the same title was published in Advances in Neural Information Processing Systems (NIPS), 2017.}
In most real word problems, parameters are uncertain at the optimization phase and decisions need to be made in the face of uncertainty. Stochastic and robust optimization are two widely used paradigms to handle uncertainty. In the stochastic optimization approach, uncertainty is modeled as a probability distribution and the goal is to optimize an expected objective~\cite{Dantzig55}. We refer the reader to Kall and Wallace~\cite{KW94}, Prekopa~\cite{Prekopa95}, Shapiro~\cite{Shapiro08}, Shapiro et al.~\cite{SDR09} for a detailed discussion on stochastic optimization. On the other hand, in the robust optimization approach, we consider an adversarial model of uncertainty using an uncertainty set and the goal is to optimize over the worst-case realization from the uncertainty set. This approach was first introduced by Soyster~\cite{SA73} and has been extensively studied in recent years. We refer the reader to Ben-Tal and Nemirovski~\cite{BN98,BN99,Ben-Tal02}, El Ghaoui and Lebret~\cite{EL97}, Bertsimas and Sim~\cite{BS03,BS04}, Goldfarb and Iyengar~\cite{GI03}, Bertsimas et al.~\cite{BBC08} and Ben-Tal et al.~\cite{BNE10} for a detailed discussion of robust optimization. In both these paradigms, computing an optimal dynamic solution is intractable in general due to the ``curse of dimensionality''.

This intractability of computing the optimal adjustable solution necessitates considering approximate solution policies. Commonly used approximations include functional policies such as static and affine policies where the decision in any period $t$ is restricted to a static or a linear function of the sample path until period $t$. Both static and affine policies have been studied extensively in the literature and can be computed efficiently for a large class of problems. While the worst-case performance of such approximate policies can be significantly bad as compared to the optimal dynamic solution, the empirical performance, especially of affine policies, has been observed to be near-optimal in a broad range of computational experiments. Our goal in this paper is to address this stark contrast between the worst-case performance bounds and near-optimal empirical performance of affine policies. 

%While it is a reasonable approach in certain settings, it is intractable in general and suffers from the ``curse of dimensionality''. Moreover, in many applications, we may not have sufficient historical data to estimate a joint probability distribution over the uncertain parameters.

 %Robust optimization leads to a tractable approach where an optimal static solution can be computed efficiently for a large class of problems. Moreover, designing an uncertainty set is significantly less challenging than estimating a joint probability distribution for high-dimensional uncertainty. However, computing an optimal adjustable (or dynamic) solution for a multi-stage problem is generally hard even in the robust optimization framework necessitating us to consider approximate policies such as static, affine and piecewise affine policies.

Specifically, we consider the following two-stage adjustable robust linear optimization problems with uncertain demand requirements:
\begin{equation}\label{eq:ar}
\begin{aligned}
z_{\sf AR}\left( \mb c, \mb d, \mb A, \mb B, {\cal U} \right) = \min_{\mb x} \; & \mb{c}^T \mb{x} + \max_{\mb{h}\in {\cal U}} \min_{\mb{y}(\mb{h})} \mb{d}^T \mb{y}(\mb{h}) \\
& \mb{A}\mb{x} + \mb{B}\mb{y}(\mb{h}) \; \geq \; \mb{h},   \; \; \; \forall \mb h \in {\cal U} \\
& \mb{y}(\mb{h}) \in  {\mathbb R}^{n}_+, \;  \; \; \forall \mb h \in {\cal U} \\
& \mb{x}  \in  {\mathbb R}^{n}_+
\end{aligned}
\end{equation}
where $\mb{A} \in {\mathbb R}_+^{m {\times} n}, \mb{c}\in {\mathbb R}^{n}_+, \mb{d}\in {\mathbb R}^{n}_+, \mb{B} \in {\mathbb R}^{m {\times} n}_+$. The right-hand-side $\mb{h}$ belongs to a compact convex uncertainty set ${\cal U}\subseteq{\mathbb R}^m_+$. The goal in this problem is to select the first-stage decision $\mb{x}$, and the second-stage recourse decision,  $\mb y(\mb{h})$, as a function of the uncertain right hand side realization, $\mb{h}$ such that the worst-case cost over all realizations of $\mb h \in {\cal U}$ is minimized. We assume without loss of generality that $\mb c= \mb e$ and $ \mb d = \bar d \cdot\mb e$ (by appropriately scaling $\mb A$ and $\mb B$). Here, $\bar d$ can interpreted as the inflation factor for costs in the second-stage.

This model captures many important applications including set cover, capacity planning and network design problems under uncertain demand. Here the right hand side, $\mb h$ models the uncertain demand and the covering constraints capture the requirement of satisfying the uncertain demand. However, the adjustable robust optimization problem \eqref{eq:ar} is intractable in general. In fact, Feige et al.~\cite{FJMM07} show that $\Pi_{\sf AR}({\cal U})$~\eqref{eq:ar} is hard to approximate within any factor that is better than $\Omega(\log n/\log \log n)$ under reasonable complexity assumptions.

Both static and affine policy approximations have been studied in the literature for~\eqref{eq:ar}. In a static solution, we compute a single optimal solution $(\mb{x},\mb{y})$ that is feasible for all realizations of the uncertain right hand side. Bertsimas et al.~\cite{BGS10}  relate the performance of static solution to the symmetry of the uncertainty set and show that it provides a good approximation to the adjustable problem if the uncertainty is close to being centrally symmetric. However, the performance of static solutions can be arbitrarily large for a general convex uncertainty set with the worst case performance being $\Omega (m)$. El Housni and Goyal \cite{elhousni2015piecewise} consider piecewise static policies for two-stage adjustable robust problem with uncertain constraint coefficients. These are a generalization of static policies where we divide the uncertainty set into several pieces and specify a static solution for each piece. However, they show that, in general, there is no piecewise static policy with a polynomial number of pieces that has a significantly better performance than an optimal static policy.

An affine policy restricts the second-stage decisions, $\mb{y}(\mb{h})$ to being an affine function of the uncertain right-hand-side $\mb{h}$, i.e.,  $\mb{y}(\mb{h})=\mb{P}\mb{h}+\mb{q}$ for some $\mb{P}\in{\mathbb R}^{n\times m}$ and $\mb{q}\in{\mathbb R}^m$ are decision variables. Affine policies in this context were introduced in Ben-Tal et al.~\cite{Ben-Tal04} and can be formulated as:
\begin{equation}\label{eq:aff}
\begin{aligned}
z_{\sf Aff}\left( \mb c, \mb d, \mb A, \mb B, {\cal U} \right) = \min_{\mb x,\mb P, \mb q} \; & \mb{c}^T \mb{x} + \max_{\mb{h}\in {\cal U}}  \mb{d}^T\left( \mb{P}\mb{h}+\mb{q}\right) \\
& \mb{A}\mb{x} + \mb{B}\left( \mb{P}\mb{h}+\mb{q}\right) \; \geq \; \mb{h},   \; \; \; \forall \mb h \in {\cal U} \\
& \mb{P}\mb{h}+\mb{q} \; \geq \; \mb{0},   \; \; \; \forall \mb h \in {\cal U} \\
& \mb{x}   \in  {\mathbb R}^{n}_+.\\
\end{aligned}
\end{equation}
An optimal affine policy can be computed efficiently for a large class of problems. Bertsimas and Goyal~\cite{BG10} show that affine policies give a $O(\sqrt{m})$-approximation to the optimal dynamic solution for~\eqref{eq:ar}. Furthermore, they show that the approximation bound $O(\sqrt m)$ is tight. However, the observed empirical performance for affine policies is near-optimal for a large set of synthetic instances of~\eqref{eq:ar}.

%The worst case scenario of problem \eqref{eq:ar} occurs on extreme points of ${\cal U}$. Therefore, given an explicit list of the extreme points of the uncertainty set ${\cal U}$, the adjustable robust optimization problem~\eqref{eq:ar} can be solved efficiently by including the second-stage decisions and the covering constraints only for the extreme points of ${\cal U}$. However, ; for example, when the number of extreme points is large or due to other structural complexities of ${\cal U}$.  

%This motivates us to consider approximate solution policies such as static policies, affine policies and piecewise policies. 

%Ben-Tal et al. \cite{bentractable} introduce a new framework to compute efficiently piecewise affine policies for the adjustable problem \eqref{eq:ar} with improved approximation bounds than affine policies for a large class of uncertainty sets.

%Note that affine policies have very good empirical performance. The usual approach is to test this policies on random instances. Our goal in this paper is to theoretically  characterize the performance of affine on random instances for a large class of distributions.

\subsection{Our Contributions}

Our goal in this paper is to address the stark contrast between the worst-case and empirical performance of affine policies and provide a more fine-grained analysis for affine policies beyond worst-case. In particular, we present a theoretical analysis of the performance of affine policies for synthetic instances of the problem generated from a probabilistic model. More specifically, we consider random instances of the two-stage adjustable problem \eqref{eq:ar} where the coefficients of the constraint matrix $\mb B$ are randomly generated and analyze the performance of affine policies for a large class of distributions. Our main contributions are summarized below.

\vspace{2mm}
\noindent {\bf Random Constraint Coefficients}. We consider probabilistic instances of~\eqref{eq:ar} where the columns of  $ \mb B$ are generated from independent multivariate distributions, (i.e., for all $j \in [n]$, column $\mb B_j$ is generated from the multivariate distribution ${\cal F}_j$  independent from the other columns) and show that affine policy is provably a good approximation with high probability with a bound that is significantly better than the worst-case bound for a large class of distributions including distributions with bounded support and distributions with gaussian and sub-gaussian tails. 

\begin{enumerate}
\item {\bf Distributions with Bounded Support}. We first consider the case where the support of distributions ${\cal F}_j$, $j \in [n]$ is bounded in, say $[0,b]^m$. For all $i \in [m]$, let ${\cal F}_{ij}$ denote the marginal distribution of $B_{ij}$ where the column $\mb B_j$ is distributed according to ${\cal F}_j$, and let  $\mu_{ij} = {\mathbb E}[ B_{ij}]$. We show that for sufficiently large values of $m$ and $n$, affine policy gives a $b/\mu$-approximation to the adjustable problem~\eqref{eq:ar} where 
\[\mu = \min_{i \in [m]} \frac{1}{n} \sum_{j=1}^n \mu_{ij}.\]   
More specifically, with probability at least $(1-1/m)$, we have that
\[ z_{\sf Aff}(\mb c, \mb d, \mb A, \mb B, {\cal U}) \leq \frac{b}{\mu (1- \epsilon)} \cdot z_{\sf AR}(\mb c, \mb d, \mb A, \mb B, {\cal U}),\]
where $\epsilon = b/\mu \sqrt{\log m/n}$ (Theorem~\ref{thm:bounded}).  

This bound is significantly better than the worst-case approximation bound of $O(\sqrt m)$ for many distributions. As an example, consider the special case where all coefficients $B_{ij}$ are i.i.d. according to some distribution with bounded support $[0,b]$ and expectation $\mu$. Then affine policy gives $b/\mu$-approximation to the two-stage adjustable problem~\eqref{eq:ar} with high probability. Moreover, if the distribution is {\em symmetric} (such as uniform or Bernoulli distribution with parameter $p=1/2$), affine policy gives a $2$-approximation for the adjustable problem~\eqref{eq:ar}.

 \item {\bf Distributions with Sub-Gaussian tails}. While the above analysis leads to a good approximation for many distributions, the ratio $\frac{b}{\mu}$ can be significantly large in general; for instance, for distributions where extreme values of the support are extremely rare and significantly far from the mean. In such instances, the bound $b/\mu$ can be quite loose. We can tighten the analysis by using the concentration properties of distributions and can extend the analysis even for the case of unbounded support. In particular, we consider the case where for all $j \in [n]$, column $\mb B_j$ is distributed according to a multivariate distribution, ${\cal F}_j$ with (possibly) unbounded support and a sub-gaussian tail independent of other columns. Then for sufficiently large values of $m$ and $n$, with probability at least $(1-1/m)$, 
\[ z_{\sf Aff}(\mb c, \mb d, \mb A, \mb B, {\cal U}) \leq O(\sqrt{\log m + \log n}) \cdot z_{\sf AR}(\mb c, \mb d, \mb A, \mb B, {\cal U}).\]
Here we assume that the parameters of the distributions are independent of the problem dimension. We prove the case of distributions with sub-gaussian tails in Theorem~\ref{thm:gaus}.
\end{enumerate}

We would like to note that the above performance bounds are in stark contrast with the worst case performance bound $O(\sqrt{m})$ for affine policies that is tight. For the random instances where columns of $\mb B$ are independent according to above distributions, with high probability the performance bound is significantly better. Therefore, our results provide a theoretical justification of the good empirical performance of affine policies and close the gap between worst case bound of $O(\sqrt{m})$ and observed empirical performance. Furthermore, surprisingly these performance bounds are independent of the structure of the uncertainty set, ${\cal U}$ unlike in previous work where the performance bounds depend on the geometric properties of ${\cal U}$. Our analysis is based on a {\em dual-reformulation} of~\eqref{eq:ar} introduced in~\cite{bertsimas2016duality} where~\eqref{eq:ar} is reformulated as an alternate two-stage adjustable optimization and the uncertainty set in the alternate formulation depends on the constraint matrix $\mb B$. Using the probabilistic structure of $\mb B$, we show that the alternate {\em dual} uncertainty set is close to a simplex for which affine policies are optimal.

We would also like to note that our performance bounds are not necessarily tight and the actual performance on particular instances can be even better. We test the empirical performance of affine policies for random instances generated according to uniform and folded normal distributions and observe that affine policies are nearly optimal with a worst optimality gap of $4\%$ (i.e. approximation ratio of $1.04$) on our test instances  as compared to the optimal adjustable solution that we compute using a mixed integer program (MIP). %for the class of Gaussian and Uniform distributions that affine policy is near optimal for the adjustable problem \eqref{eq:ar}.

%from a random  bounded distribution in $[0,b]$ with expectation $\mu$. We show that with high probability, for large values of $m$ and $n$, affine policy gives a $b/\mu$-approximation to the adjustable problem \eqref{eq:ar}. For example, Affine policy gives 2-Approximation for the adjustable problem under any uniform distribution.

%We consider a large class of distributions for $\mb B$ and give theoretical approximation bounds on the performance of affine policies. We also show empirically that affine policy is near optimal for a large class of random constraint matrix $ \mb B$ independently of the uncertainty set. We summarize our contributions below.
\vspace{2mm}
\noindent {\bf Worst-case distribution for Affine policies.} While affine policies give with high probability a good approximation  for random instances according to  a large class of commonly used distributions, we present a distribution where the performance of affine policies is $\Omega(\sqrt m)$ with high probability for instances generated from this distribution. In particular, there is no smoothed analysis for affine policies. Moreover, this bound matches the worst-case deterministic bound for affine policies. We would like to remark that in the worst-case distribution, the coefficients $B_{ij}$  depend on the dimension of the problem. This suggests that to obtain bad instances for affine policies, we need to generate instances using a structured distribution where the structure of the distribution might depend on the problem structure.

\noindent
{\bf Outline.} The rest of the paper is organized as follows. In Sect. 2, we present our results on the performance of affine policies for random instances and show that affine policies give with high probability {\em good} approximation  to $\eqref{eq:ar}$ for a large class of distributions.  In Sect. 3, we present a class of distributions and bad instances where affine policies perform poorly and match the worst-case deterministic bound.  Finally we present a
computational study to test affine policies on random instances in Sect. 4.

\section{Random instances of two-stage robust optimization problems}
In this section, we theoretically characterize the performance of affine policies for random instances of $\eqref{eq:ar}$. In particular, we consider the two-stage problem where  coefficients of constraint matrix $\mb B$ are random. 
%We would like to note that our results in this section extend as well to conic uncertainty sets such as ellipsoids (see Remark).
Our analysis of the performance of affine policies does not depend on the structure of first stage constraint matrix $\mb A$ or cost $ \mb c$.   The second-stage cost, as already mentioned, is wlog of the form $\mb d= \bar d \mb e$. Therefore, we restrict our attention only to the distribution of coefficients of the second stage matrix $\mb B$. We will use the notation $\tilde{\mb B}$ to emphasize that $\mb B$ is random. For simplicity, we refer to  $z_{\sf AR}\left( \mb c, \mb d, \mb A, \mb B, {\cal U} \right)$ as $ z_{\sf AR}\left(  \mb B \right)$ and to $z_{\sf Aff}\left( \mb c, \mb d, \mb A, \mb B, {\cal U} \right)$ as $z_{\sf Aff}\left(  \mb B \right)$.

%\subsection{Distributions with bounded support}

We first consider the case when the columns of $\tilde{\mb B}$, namely $\tilde{\mb B_j}$ for $j \in [n]$, is distributed according to a multivariate distribution ${\cal F}_j$ with bounded support in $[0,b]^m$ (for some constant $b$), independent of other columns. We compare the performance of affine policies with respect to the optimal dynamic solution and present an approximation bound that depends only on the distribution of $\tilde{\mb B}$ and holds for any uncertainty set $\cal U$.  In particular, we have the following theorem.

\begin{theorem} [Distributions with bounded support]\label{thm:bounded}
Consider the two-stage adjustable problem \eqref{eq:ar} where  $\tilde{\mb {B}_{j}}$ for $ j \in [n]$ is distributed according to a multivariate distribution, ${\cal F}_j$ with  bounded  support in $  [0,b]^m $ (for some constant $b$), independent of other columns. Let $\mathbb{E} [{\tilde{B}}_{ij} ]=\mu_{ij}$ $\forall i \in [m] \; \forall j \in [n]$.  For $n$ and $m$ sufficiently large, we have  with  probability at least $1- \frac{1}{m}$,

$$z_{\sf AR}(\tilde{\mb B}) \leq z_{\sf Aff}(\tilde{\mb B}) \leq  \frac{b}{\mu ( 1- \epsilon)} \cdot z_{\sf AR}(\tilde{\mb B})$$
where $ \mu = \min_{i \in [m]} \frac{1}{n} \sum_{j=1}^n \mu_{ij}$ and $\epsilon = \frac{b}{\mu} \sqrt{\frac{\log m}{n}}  $.
\end{theorem}

For the special case where  $\tilde{B}_{ij}$  are i.i.d. according to a bounded  distribution with support in $  [0,b] $. We have the following corollary.
\begin{corollary} \label{clr:bounded}
Consider the two-stage adjustable problem \eqref{eq:ar} where $\tilde{B}_{ij}$  are i.i.d. according to a bounded  distribution with support in $  [0,b] $ and expectation $\mu$.  For $n$ and $m$ sufficiently large, we have  with  probability at least $1- \frac{1}{m}$,

$$z_{\sf AR}(\tilde{\mb B}) \leq z_{\sf Aff}(\tilde{\mb B}) \leq  \frac{b}{\mu ( 1- \epsilon)} \cdot z_{\sf AR}(\tilde{\mb B})$$
where $\epsilon = \frac{b}{\mu} \sqrt{\frac{\log m}{n}}  $.
\end{corollary}
The above theorem and corollary show that for sufficiently large values of $m$ and $n$, the performance of affine policies is at most $b/\mu$ times the performance of an optimal adjustable solution. This shows that affine policies give a good approximation (and significantly better than the worst-case bound of $O(\sqrt m)$)  for many important distributions. We present some examples below.

\vspace{2mm}
\noindent
{\bf Example 1. [Uniform distribution]} \label{ex:uniform}
Suppose for all $ i \in [m]$ and $ j \in [n]$ $\tilde{B}_{ij}$  are i.i.d. uniform in $ [0,1] $. Then $\mu = 1/2$ and from Corollary \ref{clr:bounded} we have with  probability at least $1- 1/m$,
$$z_{\sf AR}(\tilde{\mb B}) \leq z_{\sf Aff}(\tilde{\mb B}) \leq \frac{2}{1- \epsilon} \cdot z_{\sf AR}(\tilde{\mb B})$$
where $\epsilon =2\sqrt{\log m/n}$. Therefore, for sufficiently large values of $n$ and $m$ affine policy gives a $2$-approximation to the adjustable problem in this case. Note that the approximation bound of $2$ is a conservative bound and the empirical performance is significantly better. We demonstrate this in our numerical experiments. 

\vspace{2mm}
\noindent
{\bf Example 2. [Bernoulli distribution]} \label{exr:ber}
Suppose for all $ i \in [m]$ and $ j \in [n]$, $\tilde{B}_{ij}$  are i.i.d. according to a Bernoulli distribution of parameter $p$.  Then $\mu = p$, $b=1$ and from Corollary \ref{clr:bounded} we have with  probability at least $1- \frac{1}{m}$,
$$z_{\sf AR}(\tilde{\mb B}) \leq z_{\sf Aff}(\tilde{\mb B}) \leq \frac{1}{p(1- \epsilon)}  \cdot z_{\sf AR}(\tilde{\mb B})$$
where $\epsilon = \frac{1}{p}\sqrt{ \frac{\log m}{n} }$. Therefore for constant $p$, affine policy gives a constant approximation to the adjustable problem (for example $2$-approximation for $p=1/2$).

Note that these performance bounds are in stark contrast with the worst case performance bound $O(\sqrt{m})$ for affine policies which is tight. For these random instances, the performance is significantly better. We would like to note that the above distributions are very commonly used to generate instances for testing the performance of affine policies and exhibit good empirical performance. Here, we give a theoretical justification of the good empirical performance of affine policies on such instances, thereby closing the gap between worst case bound of $O(\sqrt{m})$ and observed empirical performance.

While the approximation bound in Theorem~\ref{thm:bounded} leads to a good approximation for many distributions, the ratio $b/\mu$ can be significantly large in general. We can tighten the analysis by using the concentration properties of distributions and can extend the analysis for the case of distributions with sub-gaussian tails. In particular, we consider the case where $\tilde{\mb B}_{j}$ is generated according to a distribution with sub-gaussian tails and show a logarithmic approximation bound for affine policies. Note that we assume that the parameters of the distribution are independent of the problem dimensions. 
We have the following theorem.
\begin{theorem}[Distributions with sub-gaussian tails]
 \label{thm:gaus}
Suppose $\forall j \in [n]$, $\tilde{\mb B}_{j} = \vert \tilde{ \mb G}_{j} \vert $  such that $\tilde{\mb G}_{j} $ is a sub-Gaussian, independent of $\tilde{\mb G}_{i}$, for all $i \neq j$.   For $n$ and $m$ sufficiently large, we have  with  probability at least $1- \frac{1}{m}$,
$$z_{\sf AR}(\tilde{\mb B}) \leq z_{\sf Aff}(\tilde{\mb B}) \leq  \kappa \cdot z_{\sf AR}(\tilde{\mb B})$$
where $\kappa = O\left(  \sqrt{ \log m + \log n}  \right)$.
\end{theorem}

We would like to note that the bound of $O\left(  \sqrt{ \log m + \log n}  \right)$ depends on the dimension of the problem unlike the case of uniform bounded distributions. But, it is significantly better than the worst-case of $O(\sqrt{m})$ \cite{BG10}  for general instances. Furthermore, this bound holds for all uncertainty sets.  We would like to note though that the bounds are not necessarily tight. In fact, in our numerical experiments where the uncertainty set is a {\em budget of uncertainty},  we observe that the performance is much better than the bounds.
We discuss the intuition and the proofs of Theorem \ref{thm:bounded} and Theorem \ref{thm:gaus}  in the following subsections.

\subsection{Preliminaries}
In order to prove Theorem \ref{thm:bounded} and Theorem \ref{thm:gaus}, we need to introduce certain preliminary results. First, to develop intuition, let us consider the case of polyhedral uncertainty set ${\cal U}$, i.e.,
\begin{equation} \label{def:U}
{\cal U}=\{\mb{h}\in{\mathbb R}^m_+\;|\;    \mb R \mb h \leq \mb r \}
\end{equation}
where $\mb{R} \in {\mathbb R}^{L {\times} m}$ and $ \mb{r}\in {\mathbb R}^{L}$. This is a fairly general class of uncertainty sets that includes many commonly used sets such as {\em box uncertainty} sets and {\em budget of uncertainty} sets. In section \ref{exten}, we sketch the extension of our results to general convex uncertainty sets such as ellipsoids.
 
%In order to prove Theorem \ref{thm:bounded}, we need to introduce certain preliminary results. 
We first introduce the following formulation for the adjustable problem \eqref{eq:ar} based on ideas in Bertsimas and de Ruiter \cite{bertsimas2016duality}. 
\begin{equation}\label{eq:dual}
\begin{aligned}
z_{\sf d-AR}(\mb B)= \min_{\mb x} \; & \mb{c}^T \mb{x} + \max_{\mb{w}\in {\cal W}}    \min_{\mb{\lambda}(\mb{w})}      - (\mb A \mb x)^T \mb w   + \mb{r}^T \mb{\lambda}(\mb{w}) \\
& \mb{R}^T \mb{\lambda}(\mb{w}) \; \geq \; \mb{w},   \; \; \; \forall \mb w \in {\cal W} \\
& \mb{\lambda}(\mb{w})   \in {\mathbb R}^{L}_+, \; \forall \mb w \in {\cal W} \\
& \mb{x}  \in  {\mathbb R}^{n}_+
\end{aligned}
\end{equation}
where the set $\cal W$ is   defined as
\begin{equation} \label{def:W}
 {\cal W}=\{\mb{w}\in{\mathbb R}^m_+\;|\;    \mb B^T \mb w \leq \mb d \}.
\end{equation}
We show that the above problem is an equivalent formulation of \eqref{eq:ar}.

\begin{lemma} \label{lem:reform}
Let $z_{\sf AR}(\mb B)$ be as defined in \eqref{eq:ar} and $z_{\sf d-AR}(\mb B) $ as defined in \eqref{eq:dual}. 
Then,
$$z_{\sf AR}(\mb B)= z_{\sf d-AR}(\mb B).$$
\end{lemma}
The proof follows from  \cite{bertsimas2016duality}. For completeness, we present it in Appendix \ref{apx-proofs:lem:reform}. Reformulation \eqref{eq:dual} can be interpreted as a new two-stage adjustable problem over {\em dualized} uncertainty set ${\cal W}$ and decision $\mb{\lambda}(\mb{w})$. Following \cite{bertsimas2016duality}, we refer to \eqref{eq:dual} as the {\em dualized} formulation and to \eqref{eq:ar} as the  {\em primal} formulation.
Bertsimas and de Ruiter \cite{bertsimas2016duality} show that even the affine approximations of \eqref{eq:ar} and \eqref{eq:dual} (where recourse decisions are restricted to be affine functions of respective uncertainties) are equivalent. 
In particular, we have the following Lemma which is a restatement of Theorem 2 in \cite{bertsimas2016duality}.

\begin{lemma}{\bf (Theorem 2 in Bertsimas and de Ruiter~\cite{bertsimas2016duality})}\label{lem:berti-de Ruiter}
Let $z_{\sf d-Aff}(\mb B)$ be the objective value when $\mb{\lambda}(\mb{w}) $ is restricted to be affine function of $\mb{w}$  and $z_{\sf Aff}(\mb B) $ as defined in \eqref{eq:aff}. Then, $$z_{\sf d-Aff}(\mb B) =z_{\sf Aff}(\mb B) .$$
\end{lemma}

Bertsimas and Goyal \cite{BG10} show that affine policy is optimal for the adjustable problem \eqref{eq:ar} when the uncertainty set ${\cal U}$ is a simplex. In fact, optimality of affine policies for simplex uncertainty sets holds for more general formulation than considered in \cite{BG10}. In particular, we have the following lemma.
\begin{lemma}\label{lem:berti-goyal}
Suppose the set $ {\cal W} $ is a simplex, i.e. a convex combination of $m+1$ affinley independent points, then affine policy is optimal for the adjustable problem \eqref{eq:dual}, i.e.   $z_{\sf d-Aff}(\mb B) =z_{\sf d-AR}(\mb B) $.
\end{lemma}

The proof proceeds along similar lines as in \cite{BG10}. For completeness, we provide it in Appendix \ref{apx-proofs:lem:berti-goyal}.
 In fact, if the uncertainty set is not simplex but can be approximated by a simplex within a small scaling factor, affine policies can still be shown to be a good approximation, in particular we have the following lemma.

\begin{lemma}  \label{lem:inclusion}
Denote ${\cal W}$ the dualized uncertainty set as defined in \eqref{def:W} and suppose there exists a simplex ${\cal S}$ and $ \kappa \geq 1$ such that
$ {\cal S} \subseteq {\cal W} \subseteq \kappa\cdot {\cal S}$. Therefore,
$$ z_{\sf d-AR}(\mb B) \leq  z_{\sf d-Aff}(\mb B)  \leq \kappa \cdot z_{\sf d-AR}(\mb B).$$ 
Furthermore,
$$ z_{\sf AR}(\mb B) \leq  z_{\sf Aff}(\mb B)  \leq \kappa \cdot z_{\sf AR}(\mb B).$$
\end{lemma}
The proof of Lemma \ref{lem:inclusion} is presented in Appendix \ref{apx-proofs:lem:inclusion}.

\subsection{Proof of Theorem \ref{thm:bounded}} \label{sec:bounded}
We consider instances of problem \eqref{eq:ar} where the columns $\tilde{\mb B }_j$ are independently generated according to  bounded  distributions with support in $  [0,b]^m $. Let $\mathbb{E} [{\tilde{B}}_{ij} ]=\mu_{ij}$ for all $i\in[m], j \in [n]$ and $$ \mu = \min _{i \in [m]} \frac{1}{n} \sum_{j=1}^n \mu_{ij}.$$
Denote the dualized uncertainty set $$\tilde{{\cal W}}= \left\{\mb{w}\in{\mathbb R}^m_+\;|\;    \mb {\tilde{B}}^T \mb w \leq \bar d \cdot \mb e \right\}.$$ Our performance bound is based on showing that $\tilde{{\cal W}}$ can be sandwiched between two simplicies with a small scaling factor. In particular, consider the following simplex, 
\begin{equation}\label{def:simplex}
{\cal S}= \left\{  \mb w \in \mathbb{R}_+^m  \; \Bigg\vert \; \sum_{i=1}^m w_i \leq \frac{\bar d}{b}    \right\}.
\end{equation}
We will show that $$   {\cal S} \subseteq \tilde{{\cal W}}  \subseteq    \frac{b}{\mu ( 1- \epsilon)} \cdot {\cal S}$$ with probability at least $1 - \frac{1}{m}$
where  $\epsilon = \frac{b}{\mu} \sqrt{\frac{\log m}{n}} $. First, we show that $   {\cal S} \subseteq \tilde{{\cal W}} $. Consider any $ \mb w \in {\cal S}$. For  $j=1, \ldots , n $, we have
$$     \sum_{i=1}^m      {\tilde{B}}_{ij} w_i    \leq b  \sum_{i=1}^m w_i  \leq \bar d.$$
The first inequality holds because all components of $\mb {\tilde{B}}$ are upper bounded by $b$ and the second one follows from $ \mb w \in {\cal S}$. Hence, we have  $\mb {\tilde{B}}^T \mb w \leq \bar d \mb e$ and consequently ${\cal S} \subseteq  \tilde{{\cal W}}$.

Now, we show that the other inclusion holds with high probability. Consider any $ \mb w \in \tilde{{\cal W}}$.  We have $\mb {\tilde{B}}^T \mb w \leq \bar d \cdot \mb e $. Summing up all the inequalities and dividing by $n$, we get
\begin{equation}\label{pr:sum}
 \sum_{i=1}^m \left(  \frac{\sum_{j=1}^n \tilde{B}_{ij}}{n	} \right)  \cdot w_i \leq \bar d .
\end{equation}
The columnds of $\mb B$ are independent, hence using Hoeffding's inequality  \cite{hoeffding1963probability}  with $\tau = b\sqrt{ \frac{\log m}{n}}$ (see Appendix  \ref{apx-proofs:Hof-ineq}), we have for all $i \in [m]$,
\begin{equation*}
\mathbb{P}   \left(  \frac{\sum_{j=1}^n \tilde{B}_{ij}}{n	} - \mu_i  \geq  - \tau \right) \geq  1 - \exp \left( \frac{-2n{\tau}^2}{b^2} \right)   = 1- \frac{1}{m^2} \\
\end{equation*}
where $\mu_i = \frac{1}{n} \sum_{j=1}^n \mu_{ij}$. Then, a union bound over $i=1,\ldots,m$ gives us
\begin{equation*} 
\mathbb{P}   \left(        \frac{\sum_{j=1}^n \tilde{B}_{ij}}{n	}  \geq \mu_i - \tau  \; \; \forall i \in [m] \right)  \geq 
1- \sum_{i=1}^m \mathbb{P}   \left(        \frac{\sum_{j=1}^n \tilde{B}_{ij}}{n	}  < \mu_i - \tau   \right) \geq 1 - \sum_{i=1}^m \frac{1}{m^2}
= 1- \frac{1}{m}  .    
\end{equation*}
Therefore, with probability at least $1-\frac{1}{m} $, we have 
$$  \sum_{i=1}^m w_i    \leq   \sum_{i=1}^m      \frac{1}{\mu_i - \tau }\left(  \frac{\sum_{j=1}^n \tilde{B}_{ij}}{n	} \right)  \cdot w_i \leq \frac{ 1}{   \min_{i \in [m]} \mu_i - \tau } \cdot \sum_{i=1}^m   \left(  \frac{\sum_{j=1}^n \tilde{B}_{ij}}{n	} \right)  \cdot w_i
 \leq \frac{\bar d}{   \mu - \tau } = \frac{b}{\mu(1 - \epsilon)} \cdot \frac{\bar{d}}{b} $$
where the last inequality follows from \eqref{pr:sum}. Note that for $m$ sufficiently large , we have $ \mu - \tau> 0$. 
Then,  $ \mb w \in \frac{b}{\mu(1 - \epsilon)} \cdot {\cal S} $ for any $\mb w \in \tilde{\cal W}$. Consequently with probability at least $ 1 - 1/m$, we have $$ {\cal S}   \subseteq     \tilde{\cal W} \subseteq \frac{b}{\mu(1 - \epsilon)}  \cdot {\cal S} .$$ Finally, we apply the result of  Lemma \ref{lem:inclusion} to conclude. $ \hfill \square$

\subsection{Proof of Theorem \ref{thm:gaus}} \label{sec:unbounded}

%While the approximation bound in Theorem~\ref{thm:bounded} leads to a good approximation for many distributions, the ratio $b/\mu$ can be significantly large in general. We can tighten the analysis by using the concentration properties of distributions and can extend the analysis even for the case of distributions with unbounded support and sub-gaussian tails. In particular, we consider the case where $\tilde{\mb B}_{j}$ are distributed according to multivariate sub-gaussian and show a logarithmic approximation bound for affine policies. Note that we assume that the parameters of the multivariate Gaussian are independent of the problem dimensions. 
%We have the following theorem.
%\begin{theorem} \label{thm:gaus}
%Suppose $\forall j \in [n]$, $\tilde{\mb B}_{j} = \vert \tilde{ \mb G}_{j} \vert $  such that $\tilde{\mb G}_{j} $ is sub-Gaussian, independent of $\tilde{\mb G}_{i}$, for all $i \neq j$.   For $n$ and $m$ sufficiently large, we have  with  probability at least $1- \frac{1}{m}$,
%$$z_{\sf AR}(\tilde{\mb B}) \leq z_{\sf Aff}(\tilde{\mb B}) \leq  \kappa \cdot z_{\sf AR}(\tilde{\mb B})$$
%where $\kappa = O\left(  \sqrt{ \log m + \log n}  \right)$.
%\end{theorem}
%\proof{Proof.}
Consider instances of problem \eqref{eq:ar} where the columns $\tilde{\mb B }_j$ are independently generated according to distributions with sub-gaussian tails. In particular, we have  for all $i,j$, $\tilde{ B}_{ij} = \vert \tilde{  G}_{ij} \vert $ where  $\tilde{  G}_{ij}$ is a sub-Gaussian random variable.
Denote $$\tilde{{\cal W}}= \{\mb{w}\in{\mathbb R}^m_+\;|\;    \mb {\tilde{B}}^T \mb w \leq \bar d \cdot \mb e \}.$$ 
Our goal is to sandwich $\tilde{{\cal W}}$ between two simplicies and use Lemma \ref{lem:inclusion}.
Since  $\tilde{  G}_{ij}$ has a sub-gaussian tail, there exists positive constants $C$ and $v_{ij}$ such that for any $t>0$,
$$ \mathbb{P}  \left( \vert \tilde{G}_{ij}  \vert    \geq t \right) \leq C e^ { -  v_{ij} t^2}.$$
Therefore,
\begin{align*}
\mathbb{P} \left( \tilde{B}_{ij} \leq    \sqrt{    \frac{2 \log(mn)}{v_{ij}} }  \right) &= 1-  \mathbb{P} \left(\vert \tilde{G}_{ij} \vert > \sqrt{    \frac{2 \log(mn)}{v_{ij}} }   \right)\\
&  \geq 1-  C \exp\left( - 2 \log(mn) \right) = 1-  \frac{C}{(mn)^2}.
\end{align*}
Denote $$ \kappa = \max_{i,j}  \left(  \sqrt{    \frac{2 \log(mn)}{v_{ij}} }  \right). $$
We have $ \kappa = O \left( \sqrt{ \log m + \log n} \right)$ because $v_{ij}$ are positive constant independent of the dimensions $m$ and $n$ of the problem. Therefore by taking a union bound over $ i \in [m]$ and $ j \in [n]$ we get,
%( $ (1+x)^p \geq 1+p x\;  \forall  x > -1$).
$$\mathbb{P} \left( \tilde{B}_{ij} \leq  \kappa \; \; \forall i \in [m], \forall  j \in [n] \right)   \geq 1- \frac{C}{mn}.$$
Consider the following simplex $${\cal S}= \{  \mb w \in \mathbb{R}_+^m  \; \big\vert \; \sum_{i=1}^m w_i \leq \bar d    \}.$$ 
For any  $w \in{\cal S}$, we have with probability at least $1-\frac{C}{mn }$, 
$$  \sum_{i=1}^m \tilde{B}_{ij} w_i   \leq   \kappa \sum_{i=1}^m w_i \leq    \kappa \cdot \bar{d} \qquad \forall j \in [n] .  $$
Hence, with probability at least $1-\frac{C}{mn }$ we have,
$   {\cal S} \subseteq  \kappa \cdot \tilde{\cal W}$. 
Now, we want to find a simplex that includes $ \tilde{\cal W}$. We follow a similar approach to the proof of Theorem \ref{thm:bounded}. Consider any $ \mb w \in \tilde{{\cal W}}$.  We have similarly to equation \eqref{pr:sum}
\begin{equation} \label{toz}
 \sum_{i=1}^m \left(  \frac{\sum_{j=1}^n \tilde{B}_{ij}}{n	} \right)  \cdot w_i \leq \bar d .
\end{equation}
We have the following concentration inequality for non-negative random variables,
\begin{equation*}
\mathbb{P}   \left(  \frac{\sum_{j=1}^n \tilde{B}_{ij}}{n	}   \geq \mu_i  - \tau_i \right) \geq     1- \exp \left(  \frac{-n \tau_i^2}{2 \sigma_i^2}                \right)         =1- \frac{1}{m^2} \\
\end{equation*}
where $\tau_i = 2 \sigma_i \sqrt{ \frac{\log m}{n}}$,  $ \mu_i =\frac{1}{n} \sum_{j=1}^n \mathbb{E} [\tilde{B}_{ij}] $ and $ \sigma_i^2 = \max_j {\sf Var}[\tilde{B}_{ij}]   $.
Then,  a union bound over $i \in [m]$ gives 
\begin{equation*} \label{pr:prob}
\mathbb{P}   \left(        \frac{\sum_{j=1}^n \tilde{B}_{ij}}{n	}  \geq \mu_i - \tau_i,  \; \; \forall i \in [m] \right)     \geq 1- \frac{1}{m} ,    
\end{equation*}
which implies 
$$ \mathbb{P}   \left(        \frac{\sum_{j=1}^n \tilde{B}_{ij}}{n	}  \geq {\kappa}',  \; \; \forall i \in [m] \right)     \geq 1- \frac{1}{m}  .    $$
where ${\kappa}' =    \max_{ i \in [m]} ( \mu_i- \tau_i)$.
 Therefore, combining this result with inequality \eqref{toz}, we have with probability at least $1-\frac{1}{m}$,
$   \tilde{\cal W} \subseteq \frac{1}{ \kappa' } \cdot \cal S$. Denote, ${\cal S'} =  \frac{1}{\kappa} {\cal S}$. We have shown that with probability at least $1- C/mn$, $ {\cal S'} \subseteq \tilde{\cal W}$. Therefore, we have with probability at least $1-\frac{1}{m}$, $$ {\cal S'} \subseteq \tilde{\cal W} \subseteq \frac{\kappa}{\kappa'}  \cdot {\cal S'} $$ where 
$$ \frac{\kappa}{\kappa'}    = O\left(  \sqrt{ \log m + \log n }  \right),$$ for sufficiently large values of $m$ and $n$.
We finally use  Lemma \ref{lem:inclusion} to conclude.
\hfill
\Halmos
\endproof

%We can extend the analysis and show a similar bound for the class of distributions with sub-Gaussian tails. The bound in Theorem \ref{thm:gaus} also hold if the columns of $\tilde{\mb B}$ are generated from independent multivariate Gaussian and the proof is along the same line as in \ref{thm:gaus}. 

\subsection{Extension to general convex uncertainty sets} \label{exten}

In this section, we show that our results  of Theorem \ref{thm:bounded} and Theorem \ref{thm:gaus} hold as well for general convex uncertainty sets $\cal U$ including ellipsoids and norm-ball sets that are widely used in  robust optimization. This is based on approximating a convex uncertainty set by a polyhedral set (possibly given by an exponential number of  inequalities). In fact, in Section \ref{sec:bounded} and Section \ref{sec:unbounded}, we prove Theorem \ref{thm:bounded} and Theorem \ref{thm:gaus} for the case of polyhedral uncertainty set $ \cal U$. Note that the approximation bounds are independent from the description of ${\cal U}$ and depend only on the distribution of $\tilde{\mb B}$.

Now, consider a genral convex  uncertainty set ${\cal U} \subseteq \mathbb{R}^m$. For any $\epsilon > 0$, Deville et al. \cite{deville1998analytic} show that  there exists a polyhedral set ${\cal V}$  (see Theorem 1.1 in \cite{deville1998analytic})  such that 
\begin{equation} \label{inclu}
{\cal V} \subseteq {\cal U} \subseteq (1+ \epsilon) \cdot {\cal V}.
\end{equation}
Note that the number of polyhedral inequalities that describes $ {\cal V}$  could be exponential in $m$ and $1/ \epsilon$. Consider instances of the two-stage adjustable problem \eqref{eq:ar}  with random second-stage matrix $\tilde{\mb B}$. Denote $\beta$ the approximation bound given by Theorem \ref{thm:bounded} or Theorem \ref{thm:gaus} on the performance of affine policies for polyhedral uncertainty sets. Note that $\beta$ depends only on the distribution of $\tilde{\mb B}$ and does not depend on the description of the polyhedral uncertainty set. Therefore,
$$  z_{\sf Aff}(\tilde{\mb B}, {\cal V})     \leq  \beta  \cdot  z_{\sf AR}(\tilde{\mb B}, {\cal V}), $$
where we use the notation $z(\tilde{\mb B}, {\cal V})$ to denote the adjustable or affine problem with random matrix $\tilde{\mb B}$ and uncertainty set $ {\cal V}$.
Combining the above inequality with \eqref{inclu}, we get
$$z_{\sf Aff}(\tilde{\mb B}, {\cal U}) \leq  ( 1+ \epsilon) \cdot  z_{\sf Aff}(\tilde{\mb B}, {\cal V})     \leq  \beta ( 1+ \epsilon)  \cdot  z_{\sf AR}(\tilde{\mb B}, {\cal V}) \leq    \beta ( 1+ \epsilon)  \cdot  z_{\sf AR}(\tilde{\mb B}, {\cal U}).$$
Since $\epsilon >0 $ could be chosen arbitrary small, then
$$z_{\sf Aff}(\tilde{\mb B}, {\cal U})  <   \beta   \cdot  z_{\sf AR}(\tilde{\mb B}, {\cal U}).$$
i.e., the same approximation bounds of Theorem \ref{thm:bounded} and Theorem \ref{thm:gaus} hold as well for general convex uncertainty sets.
%\hfill
%\Halmos
%\endproof

\begin{remark} 
We would like to note that our results extend as well for two-stage  robust optimization problems \eqref{eq:ar} where  the constraints matrices $\mb A$ and $\tilde{\mb B}$ could possibly have some negative components. In fact, the non-negativity assumption on $\mb A$  could be relaxed without loss of generality since our analysis in the paper depends only on the second stage matrix $\tilde{\mb B}$. We can relax the non-negativity of $\tilde{\mb B}$ under two  assumptions:
\begin{enumerate}
\item The affine  problem $z_{\sf Aff}(\tilde{\mb B})$  is feasible.
\item For each row $i \in [m]$ of $\tilde{\mb B}$, $$\mu_i =\frac{1}{n} \sum_{j=1}^n \mathbb{E} [\tilde{B}_{ij}] > 0.$$
\end{enumerate}
In fact, in the proof of  Theorem \ref{thm:bounded} and Theorem \ref{thm:gaus}, we did not require  the matrix $\tilde{\mb B}$ to be non-negative but we used only the fact that $\mu_i - \tau_i \geq 0$ for  small enough $\tau_i$. Hence, our second-stage matrix $\tilde{\mb B}$ could have negative components as long as $\mu_i > 0$ for all rows $i=1,\ldots,m$.  On the other hand, Assumption 1 is required because feasibility of the affine problem is not necessary guaranteed if we relax the non-negativity  of both  
 matrices $\mb A$ and $\tilde{\mb B}$. 
\end{remark}

\section{Family of worst-case distribution}

For any $m$ sufficiently large, the authors in \cite{BG10} present an instance where affine policy is $\Omega( m^{\frac{1}{2}- \delta }    )$ away from the optimal adjustable solution. The parameters of the instance in \cite{BG10} were carefully chosen to achieve the gap $\Omega(m^{\frac{1}{2}- \delta})$. In this section, we show that the family of worst-case instances is not a measure zero set. In fact, we exhibit a distribution and an uncertainty set such that a random instance, $\tilde{\mb B}$ sampled from that distribution achieves a worst-case bound of $\Omega(\sqrt{m})$ with high probability. The coefficients $\tilde{B}_{ij}$ in our bad family of instances are independent but they depend on the dimension of the problem. The instance can be given as follows.
\begin{equation} \label{ex:bad:dist}
\begin{aligned}
&n=m, \; \;  \mb A= \mb 0, \; \;  \mb c= \mb 0, \; \;  \mb d = \mb e \\
&{\cal U} = {\sf{conv}} \left( \mb 0, \mb e_1 , \ldots ,\mb e_m, \mb \nu_1, \ldots, \mb \nu_m \right) \; 
\text{        where   } \mb \nu_i =  \frac{1}{\sqrt{m}}( \mb e- \mb e_i) \; \forall i \in [m]. \\
&\tilde{B}_{ij}= \left\{
    \begin{array}{ll}
        1 & \mbox{if} \; \; i=j \\
          \frac{1}{\sqrt{m}} \cdot \tilde{u}_{ij}		& \mbox{if}  \; \;   i \neq j \\
    \end{array}
\right.
 \text{where for all   } i \neq j, \tilde{u}_{ij} \text{   are i.i.d. uniform} [0,1].
 \end{aligned}
\end{equation}
\begin{theorem} \label{thm:bad:dist}
For the instance defined in \eqref{ex:bad:dist}, we have with probability at least $1-1/m$,   
$$z_{\sf Aff}({\tilde{\mb B}})= \Omega ( \sqrt{m}) \cdot  z_{\sf AR}({\tilde{\mb B}}).$$
\end{theorem}
As a byproduct, we also tighten the lower bound on the performance of affine policy to $\Omega(\sqrt m)$ improving from the lower bound of $\Omega( m^{\frac{1}{2}- \delta} ) $ in \cite{BG10}. We would like to note that both uncertainty set and distribution of coefficients  in our instance \eqref{ex:bad:dist} are carefully chosen to achieve the worst-case gap. Our analysis suggests that to obtain bad instances for affine policies, we need to generate instances using a structured distribution as above and it may not be easy to obtain bad instances in a completely random setting as observed in extensive empirical studies. 

%\section{Proofs of Theorem \ref{thm:bad:dist}}
%\subsection*{Proof of Theorem \ref{thm:bad:dist}}
To prove Theorem \ref{thm:bad:dist}, we introduce the following lemma which shows a deterministic bad instance where the optimal affine solution is $\Theta (\sqrt{m})$ away from the optimal adjustable solution.
\begin{lemma} \label{lem:worst-case}
Consider the two-stage adjustable problem \eqref{eq:ar} where:
\begin{equation}
\begin{aligned}
&n=m, \; \;  \mb A= \mb 0, \; \;  \mb c= \mb 0, \; \;  \mb d = \mb e \\ 
&{\cal U} = {\sf{conv}} \left( \mb 0, \mb e_1 , \ldots ,\mb e_m, \mb \nu_1, \ldots, \mb \nu_m \right) \; 
\text{        where   } \mb \nu_i =  \frac{1}{\sqrt{m}}( \mb e- \mb e_i) \; \forall i \in [m]. \\ 
&B_{ij}= \left\{
    \begin{array}{ll}
        1 & \mbox{if} \; \; i=j \\
          \frac{1}{\sqrt{m}}	& \mbox{if}  \; \;   i \neq j 
    \end{array}
\right.
 \end{aligned} \label{ex:deter}
\end{equation}
%$n=m,  \mb c = \mb 0,  \; \mb d  = \mb e, \mb A= \mb 0$,
%\begin{equation} \label{matrix:B}
 %\mb B= \left(
%\begin{matrix}
%1 &  \frac{1}{\sqrt{m}} & \ldots & \frac{1}{\sqrt{m}}\\
%\frac{1}{\sqrt{m}} &  1 & \ldots & \frac{1}{\sqrt{m}}\\
%\vdots & \vdots & \ddots & \vdots\\
%\frac{1}{\sqrt{m}}  &   \frac{1}{\sqrt{m}}       &\ldots & 1
%\end{matrix}
%\right)
%\end{equation}
%\begin{equation} \label{matrix:B}
%     {B}_{ij}= \left\{
%    \begin{array}{ll}
%        1 & \mbox{if} \; \; i=j \\
%          \frac{1}{\sqrt{m}} 		& \mbox{if}  \; \;   i \neq j \\
%    \end{array}
%\right.
%\end{equation}
%and the uncertainty set is defined as
%\begin{equation} \label{eq:worstU}
%{\cal U} = {\sf{conv}} \left( \mb 0, \mb e_1 , \ldots ,\mb e_m, \mb \nu_1, \ldots, \mb \nu_m \right)
%\end{equation}
%where $\mb \nu_i =  \frac{1}{\sqrt{m}}( \mb e- \mb e_i)$ for $i=1, \ldots,m$.
Then, $z_{\sf Aff}({\mb B})= \Omega ( \sqrt{m}) \cdot  z_{\sf AR}({\mb B}).$
\end{lemma}

\proof{Proof.}
First, let us prove that $z_{\sf AR}({\mb B}) \leq 1$. It is sufficient to define an adjustable solution only for the extreme points of ${\cal U}$ because the constraints are linear. We define the following solution for all $i \in [m]$,
$$  \mb x = \mb 0 , \qquad \mb y ( \mb 0) = \mb 0, \qquad  \mb y ( \mb e_i) = \mb e_i, \qquad  \mb y ( \mb \nu_i) = \frac{1}{m} \mb e.$$
We have $\mb B \mb y ( \mb 0) = \mb 0$.  For all $i \in [m]$,  
$$\mb B \mb y ( \mb e_i) = \mb e_i + \frac{1}{\sqrt{m}} ( \mb e - \mb e_i) \geq  \mb e_i$$
and 
$$\mb B \mb y (\mb \nu_i) = \frac{1}{m} \mb B \mb e = \left( \frac{1}{m}+ \frac{m-1}{m \sqrt{m}} \right) \mb e \geq \frac{1}{\sqrt{m}} \mb e    \geq \mb \nu_i.$$
Therefore, the solution defined above is feasible. Moreover, the cost of our feasible solution is $1$ because for all $i \in [m]$,   we have
$$ \mb d^T \mb y ( \mb e_i)= \mb d^T \mb y ( \mb \nu_i)= 1.$$ 
Hence, $z_{\sf AR}({\mb B}) \leq 1.$
Now, it is sufficient to prove that  $z_{\sf Aff}({\mb B})= \Omega ( \sqrt{m})$. From Lemma 8 in Bertsimas and Goyal \cite{BG10}, since our instance is symmetric, i.e., the uncertainty set ${\cal U}$ and the dualized uncertainty set ${\cal W}$ are permutation invariant,  there exists an optimal solution for the affine problem \eqref{eq:aff} of the following form $ \mb y( \mb h) = \mb P \mb h + \mb q$ for $\mb h \in {\cal U}$ where
\begin{equation} \label{matrix:P}
 \mb P= \left(
\begin{matrix}
\theta &  \mu & \ldots & \mu \\
\mu &  \theta & \ldots & \mu \\
\vdots & \vdots & \ddots & \vdots\\
\mu  &   \mu    &\ldots & \theta
\end{matrix}
\right)
\end{equation}
and $ \mb q = \lambda \mb e$. We have $ \mb y(\mb 0) = \lambda \mb e \geq \mb 0$ hence 
\begin{equation}\label{eq:lambda}
\lambda \geq 0.
\end{equation}
We know that 
\begin{equation}\label{eq:lem:1}
z_{\sf Aff}({\mb B}) \geq \mb d^T \mb y( \mb 0) = \lambda m.
\end{equation}

\vspace{2mm}
\noindent
{\bf Case 1:} If $\lambda \geq \frac{1}{6 \sqrt{m}}$, then from \eqref{eq:lem:1} we have $z_{\sf Aff}({\mb B}) \geq \frac{\sqrt{m}}{6}$.

\vspace{2mm}
\noindent
{\bf Case 2:} If $ \lambda \leq \frac{1}{6 \sqrt{m}}$. We have,
$$ \mb y( \mb e_1) = ( \theta+ \lambda ) \mb e_1 + ( \mu+ \lambda) ( \mb e -  \mb e_1).$$
By feasibility of the solution, we have $ \mb B \mb y ( \mb e_1) \geq \mb e_1$, hence
$$ (\theta+ \lambda) +\frac{1}{\sqrt{m}} (m-1)(\mu+ \lambda) \geq 1.$$
Therefore $ \theta + \lambda \geq \frac{1}{2}$ or $\frac{1}{\sqrt{m}} (m-1)(\mu + \lambda) \geq \frac{1}{2}$.

\vspace{2mm}
\noindent
{\bf Case 2.1:} Suppose $\frac{1}{\sqrt{m}} (m-1)(\mu + \lambda) \geq \frac{1}{2}$.
Therefore,
$$ z_{\sf Aff}({\mb B}) \geq \mb d^T \mb y (\mb e_1) = \theta+ \lambda + (m-1)(\mu+ \lambda) \geq \frac{\sqrt{m}}{2}.$$
where the last inequality holds because $\theta+ \lambda \geq 0 $ as $\mb y( \mb e_1) \geq \mb 0$.

\vspace{2mm}
\noindent
{\bf Case 2.2:} Now suppose we have the other inequality i.e., $ \theta + \lambda \geq \frac{1}{2}$. Recall that we have $ \lambda \leq \frac{1}{6\sqrt{m}}$ as well. Therefore,
$$ \theta \geq \frac{1}{2}- \frac{1}{6\sqrt{m}} \geq \frac{1}{3}.$$
We have,
$$ \mb y ( \mb \nu_1 )  =  \frac{1}{\sqrt{m}}     \left( ( \theta + (m-2) \mu ) ( \mb e- \mb e_1)   + (m-1) \mu \mb e_1 \right) + \lambda \mb e. $$
Therefore,
\begin{align} \label{eq:up} \nonumber
 z_{\sf Aff}({\mb B}) \geq \mb d^T \mb y (\mb \nu_1) &= \frac{1}{\sqrt{m}} ( (m-1) \theta +(m-1)^2 \mu) + \lambda m \\ 
 & \geq    \frac{m-1}{\sqrt{m}} \left(  \frac{1}{3} + (m-1) \mu \right). 
\end{align}
where the last inequality follows from $\lambda \geq 0$ and $ \theta \geq \frac{1}{3}.$

\vspace{2mm}
\noindent
{\bf Case 2.2.1:} If $ \mu \geq 0$ then from \eqref{eq:up}
$$ z_{\sf Aff}({\mb B}) \geq  \frac{m-1}{3\sqrt{m}} = \Omega ( \sqrt{m}).$$

\vspace{2mm}
\noindent
{\bf Case 2.2.2:} Now suppose that $ \mu < 0$, by non-negativity of $ \mb y ( \mb \nu_1) $ we have,
$$ \frac{m-1}{\sqrt{m}} \mu + \lambda \geq 0$$ 
i.e., $$ \mu \geq \frac{-\lambda \sqrt{m}}{m-1},$$
and from \eqref{eq:up}
\begin{align*}
 z_{\sf Aff}({\mb B})  &\geq    \frac{m-1}{\sqrt{m}} \left(  \frac{1}{3} + (m-1) \mu \right) \\
 & \geq  \frac{m-1}{\sqrt{m}}\left(  \frac{1}{3} - \lambda \sqrt{m}\right) \\
& \geq  \frac{m-1}{\sqrt{m}}\left(  \frac{1}{3} - \frac{1}{6}\right) =   \frac{m-1}{6 \sqrt{m}} = \Omega ( \sqrt{m}). 
\end{align*}
We conclude that in all cases $z_{\sf Aff}({\mb B})= \Omega ( \sqrt{m})$ and consequently $z_{\sf Aff}({\mb B})= \Omega ( \sqrt{m}) \cdot  z_{\sf AR}({\mb B}).$
\hfill
\Halmos
\endproof

\subsection*{Proof of Theorem \ref{thm:bad:dist}}\label{apx-proofs:thm:bad:dist}
\proof{Proof.}
Denote  $${\cal W}=\{\mb{w}\in{\mathbb R}^m_+\;|\;    \mb B^T \mb w \leq  \bar d  \mb e\}$$ and $$\tilde{\cal W}=\{\mb{w}\in{\mathbb R}^m_+\;|\;    \tilde{\mb B}^T \mb w \leq  \bar d  \mb e\}$$
 where $\mb B$  is defined in \eqref{ex:deter} and $\tilde{\mb B}$  is defined in \eqref{ex:bad:dist}. We know  for all $i,j$ in $\{1,\ldots,m\}$ that $\tilde{B}_{ij} \leq {B}_{ij}$. Hence, for any $\mb w \in {\cal W}$, we have $  \tilde{\mb B}^T \mb w \leq   {\mb B}^T \mb w   \leq \bar d \mb e$. Therefore $ \mb w \in \tilde{\cal W}$ and consequently 
 ${\cal W} \subseteq \tilde{\cal W}$. 
 Now, suppose $\mb w \in \tilde{\cal W}$, we have for all $i \in [m],$
  \begin{equation}\label{eq:thm:worst3}
 w_i + \frac{1}{\sqrt{m}} \sum_{\underset{j \neq i}{j=1} }^m \tilde{u}_{ji} w_j \leq \bar d .
 \end{equation}
 By taking the sum over $i \in [m]$, dividing by $m$ and rearranging, we get
 \begin{equation}\label{eq:thm:worst}
 \sum_{i=1}^m w_i \left(\frac{1}{m}  + \frac{1}{m\sqrt{m}} \sum_{\underset{j \neq i}{j=1} }^m \tilde{u}_{ij}   \right)   \leq \bar d.
 \end{equation}
We apply Hoeffding's inequality \cite{hoeffding1963probability} (see appendix \ref{apx-proofs:Hof-ineq}) with $\tau=\sqrt{\frac{\log m}{m-1}}$,
\begin{equation*}
\mathbb{P}   \left(  \frac{\sum_{\underset{j \neq i}{j=1} }^m \tilde{u}_{ij}}{m-1	}  \geq   \frac{1}{2}- \tau \right) \geq  1 - \exp \left( -2(m-1){\tau}^2 \right)   = 1- \frac{1}{m^2},
\end{equation*}
and we take a union bound over $j=1,\ldots,m$, we get
\begin{equation} \label{eq:thm:worst2}
\mathbb{P}   \left(   \frac{\sum_{\underset{j \neq i}{j=1} }^m \tilde{u}_{ij}}{m-1}  \geq \frac{1}{2} - \tau  \; \; \forall j=1,\ldots,m \right)  \geq \left(  1-       \frac{1}{m^2}   \right)^m     \geq 1- \frac{1}{m}, 
\end{equation}
where the last inequality follows from Bernoulli's inequality. Therefore, we conclude from \eqref{eq:thm:worst} and \eqref{eq:thm:worst2}, that with probability at least $1- \frac{1}{m}$ we have $$ \beta \sum_{i=1}^m w_i \leq \bar{d}$$ where $$\beta = \frac{1}{m}+ \frac{m-1}{m\sqrt{m}}( \frac{1}{2}- \tau) \geq \frac{1}{4\sqrt{m}}$$ for $m$ sufficiently large.
Note from \eqref{eq:thm:worst3} that for all $i$ we have $w_i \leq \bar d$. Hence with probability at least $1- \frac{1}{m}$, we have for all $i=1,\ldots,m$
$$ \mb B^T_i \mb w = w_i + \frac{1}{\sqrt{m}} \sum_{\underset{j \neq i}{j=1} }^m  w_j  \leq \bar{d}+ \frac{\bar{d}}{\beta \sqrt{m}} \leq 5 \cdot \bar{d}.$$
Therefore, $\mb w \in 5 \cdot {\cal W}$ for any $\mb w$ in ${\cal W}$ and consequently we have with probability at least $1 - \frac{1}{m}$ that, $ \tilde{\cal W} \subseteq 5 \cdot {\cal W}$. All together we have proved with probability at least $1 - \frac{1}{m},$ that $${\cal W}  \subseteq   \tilde{\cal W} \subseteq 5 \cdot {\cal W}.$$ 
This implies with probability at least $1 - \frac{1}{m}$, that $z_{\sf d-Aff}(\tilde{\mb B}) \geq z_{\sf d-Aff}({\mb B})$ and   $\   z_{\sf d-AR}({\mb B})  \geq \frac{z_{\sf d-AR}(\tilde{\mb B}) }{5}$.
We know from from Lemma \ref{lem:berti-de Ruiter} and Lemma \ref{lem:reform} that the dualized and primal are the same both for the adjustable problem and affine problem. Hence,
with probability at least $1 - \frac{1}{m}$, we have  $z_{\sf Aff}(\tilde{\mb B}) \geq z_{\sf Aff}({\mb B})$ and   $\   z_{\sf AR}({\mb B})  \geq \frac{z_{\sf AR}(\tilde{\mb B}) }{5}$.

Moreover, we know from Lemma \ref{lem:worst-case} that $z_{\sf Aff}({\mb B})  \geq  \Omega( \sqrt{m} ) \cdot   z_{\sf AR}({\mb B})$. Therefore, with probability at least $1 - \frac{1}{m}$,
$$z_{\sf Aff}(\tilde{\mb B}) \geq  \Omega( \sqrt{m} ) z_{\sf AR}(\tilde{\mb B}) .$$ 

\hfill
\Halmos
\endproof

\section{Performance of affine policy: Empirical study}
 In this section, we present a computational study to test the empirical performance of affine policy for the two-stage adjustable problem \eqref{eq:ar} on random instances. 
 
 \vspace{1mm}
\noindent {\bf Experimental setup.}  We consider two classes of distributions for generating random instances: $i)$ Coefficients of $\tilde{\mb B}$ are i.i.d. uniform $[0,1]$,  and $ii)$ Coefficients of $\tilde{\mb B}$ are absolute value of i.i.d. standard Gaussian. We consider the following {\em budget of uncertainty} set.
\begin{equation}\label{set-ones}
 {\cal U}= \left\{ \mb h \in [0,1]^m \;  \bigg\vert \;     \sum_{i=1}^m h_i  \leq \sqrt{m}        \right\}.
\end{equation}
Note that the set \eqref{set-ones} is widely used in both theory and practice and arises naturally as a consequence of concentration of sum of independent uncertain demand requirements. We would like to also note that the adjustable problem over this budget of uncertainty, ${\cal U}$ is hard to approximate within a factor better than $O(\frac{\log n}{\log \log n})$~\cite{FJMM07}. We consider
$n=m, \mb d = \mb e$. Also, we consider $\mb c= \mb 0, \mb A= \mb 0$. We restrict to this case in order to compute the optimal adjustable solution in a reasonable time by solving a single MIP. For the general problem, computing the optimal adjustable solution requires solving a sequence of MIPs each one of which is significantly challenging to solve. We would like to note though that our analysis does not depend on the first stage cost $\mb c$ and matrix $\mb A$ and affine policy can be computed efficiently even without this assumption.
We consider values of $m$ from $10$ to $ 50$ and consider $20$ instances for each value of $m$. We report the ratio $r =  z_{\sf{Aff}}(\tilde{\mb B})   / z_{\sf AR}  (\tilde{\mb B}) $ in Table \ref{tab}. In particular, for each value of $m$, we report the average ratio $r_{\sf avg}$, the maximum ratio  $r_{\sf max}$, the running time of adjustable policy $T_{\sf AR}(s)$ and the running time of affine policy $T_{\sf Aff}(s)$. We first give a compact LP formulation for the affine problem \eqref{eq:aff} and a compact MIP formulation for the separation of the adjustable  problem\eqref{eq:ar}. %Then, we compare the affine policy and the adjustable policy both in terms of objective value and computation time for different random instances. 
 
 \vspace{1mm}
 \noindent {\bf LP formulations for the affine policies.}
 The affine problem \eqref{eq:aff} can be reformulated as follows
  \iffalse%%%%%%%%%%%%%
\[
z_{\sf Aff}({\mb B}) = \min \; \left\{ \mb{c}^T \mb{x} + z \; \left | \; \begin{array}{ll}
& z \geq   \mb{d}^T\left( \mb{P}\mb{h}+\mb{q}\right)  \; \; \; \forall \mb h \in {\cal U} \\
 & \mb{A}\mb{x} + \mb{B}\left( \mb{P}\mb{h}+\mb{q}\right) \; \geq \; \mb{h}   \; \; \; \forall \mb h \in {\cal U} \\
 & \mb{P}\mb{h}+\mb{q} \; \geq \; \mb{0}   \; \; \; \forall \mb h \in {\cal U} \\
 & \mb{x}   \in  {\mathbb R}^{n}_+
 \end{array}
 \right. \right \}.
\]
\fi%%%%%%%%%%%%%%%

\begin{equation*}
\begin{aligned}
z_{\sf Aff}({\mb B}) = \min_{\mb x} \; & \mb{c}^T \mb{x} + z\\
& z \geq   \mb{d}^T\left( \mb{P}\mb{h}+\mb{q}\right)  \; \; \; \forall \mb h \in {\cal U} \\
& \mb{A}\mb{x} + \mb{B}\left( \mb{P}\mb{h}+\mb{q}\right) \; \geq \; \mb{h}   \; \; \; \forall \mb h \in {\cal U} \\
& \mb{P}\mb{h}+\mb{q} \; \geq \; \mb{0}   \; \; \; \forall \mb h \in {\cal U} \\
& \mb{x}   \in  {\mathbb R}^{n}_+.\\
\end{aligned}
\end{equation*}

Note that this formulation has infinitely many constraints but we can write a compact LP formulation using standard techniques from duality. For example, the first constraint is equivalent to 
$$z - \mb d^T \mb q \geq \max \; \{ \mb{d}^T\mb{P}\mb{h}\; | \;  \mb R  \mb h \leq \mb r, \; \mb h \geq \mb 0\}.$$ 
By taking the dual of the maximization problem, the constraint becomes $$z - \mb d^T \mb q \geq \min \;     \{\mb{r}^T\mb v \; | \; \mb R^T \mb v \geq \mb P^T \mb d, \; \mb v \geq \mb 0\}.$$
%$$ z - \mb d^T \mb q \geq \underset {\mb v \geq \mb 0 }{\underset{   \mb R^T \mb v \geq \mb P^T \mb d  }\min } \;      \mb{r}^T\mb v $$
We can then drop the min and introduce $\mb v $ as a variable, hence we obtain the following linear constraints 
$$ z- \mb d^T \mb q \geq \mb r^T \mb v, \qquad \mb R^T \mb v \geq \mb P^T \mb d, \qquad \mb{v}  \geq \mb 0.$$ We can apply the same techniques for the other constraints. The complete LP formulation and its  proof of correctness is presented in Appendix \ref{apx-LP-MIP}.

\vspace{1mm}
\noindent {\bf MIP Formulation for the adjustable problem~\eqref{eq:ar}}. 
For the adjustable problem~\eqref{eq:ar}, we show that the separation problem \eqref{eq:sep} can be formulated as a mixed integer program (MIP). The separation problem can be formulated as follows:
 \noindent
%{\bf Separation Problem.} 
Given $ \hat{\mb x}$ and $ \hat{z}$ decide whether
\begin{equation}\label{eq:sep}
\max \; \{ (\mb h - \mb A  \hat{\mb x})^T \mb w \; | \; \mb w \in {\cal W},  \mb h \in {\cal U}\}  > \hat{z} 
\end{equation}
 
 The correctness of formulation \eqref{eq:sep} follows from equation \eqref{eq:proof:sep} in the proof of Lemma \ref{lem:reform} in  Appendix  \ref{apx-proofs:lem:reform}.  The constraints in~\eqref{eq:sep} are linear but the objective function contains a bilinear term, ${\mb h}^T \mb w$. We linearize this using a standard {\em digitized reformulation}. In particular, we consider finite bit representations of continuous variables, $h_i$ nd $w_i$ to desired accuracy and introduce additional binary variables, $\alpha_{ik}$, $\beta_{ik}$ where $\alpha_{ik}$ and $\beta_{ik}$ represents the $k^{th}$ bits of $h_i$ and $w_i$ respectively. Now, for any $i\in [m]$, $h_i \cdot w_i$ can be expressed as a bilinear expression with products of binary variables, $\alpha_{ik} \cdot \beta_{ij}$ which can be linearized using additional variable $\gamma_{i jk}$ and standard linear inequalities: $ \gamma_{ij k} \leq \beta_{ij}$, $ \gamma_{ijk} \leq \alpha_{ik} $, $ \gamma_{ijk} +1 \geq \alpha_{ik} + \beta_{ij} $. The complete MIP formulation and the proof of correctness is presented in Appendix  \ref{apx-LP-MIP}.

 \iffalse%%%%%%%%%%%%%%%%%%%%%%%%%%%%%
  Let $ \mb h = \sum_{i=1}^m h_i \mb e_i$. For all $i \in [m]$ we digitize the component $ h _i$ as follows 
$$h _i= \sum_{k= -\Delta_{\cal U}}^s \frac{\alpha_{ik}}{2^k}$$
where $s= \lceil{ \log_2 \left( \frac{m}{\epsilon } \right) \rceil}$,  $\Delta_{\cal U}$ is an upper bound on any  $ h_i$ and $\alpha_{ik}$ are binary variables. This digitization gives an approximation to $h_i$ within $\frac{\epsilon }{m}$ which translates to an accuracy of $\epsilon$ in the objective function. We have
$$ \mb h = \sum_{i=1}^m \sum_{k= -\Delta_{\cal U}}^s \frac{\alpha_{ik}}{2^k} \cdot  \mb e_i \qquad 
\text{and similarly,} \qquad 
 \mb w = \sum_{i=1}^m \sum_{j= -\Delta_{\cal W}}^s \frac{\beta_{ij}}{2^j} \cdot  \mb e_i $$
where $\Delta_{\cal W}$ is an upper bound on any component of $ w \in {\cal W}$. Therefore, the first term in the objective function becomes
$$   \mb h^T \mb w=         \sum_{i=1}^m \sum_{j= -\Delta_{\cal W}}^s  \sum_{k= -\Delta_{\cal U}}^s  \frac{1}{2^{j+k}} \cdot \alpha_{ik}  \beta_{ij}$$
The final step is to linearize the term $\alpha_{ik}  \beta_{ij}$. We set, $\alpha_{ik}  \beta_{ij}=\gamma_{ijk}$ where again $\gamma_{ijk}$ is a binary variable. We can express  $\gamma_{ijk}$ using only linear constraints as follows  $ \gamma_{ijk} \leq \beta_{ij}$, $ \gamma_{ijk} \leq \alpha_{ik} $ and $ \gamma_{ijk} +1 \geq \alpha_{ik} + \beta_{ij} $. The complete MIP formulation  in presented in Appendix  \ref{apx-LP-MIP}. 
\fi%%%%%%%%%%%%%%%%%%%%%%%%%%%%%%%%%

For general $\mb A \neq 0$, we need to solve a sequence of MIPs to find the optimal adjustable solution. In order to compute the optimal adjustable solution in a reasonable time, we assume $\mb A=0, \mb c=0$ in our experimental setting so that we only need to solve one MIP.% the adjustable problem is equivalent to \eqref{eq:sep} which can be formulated as one MIP.

\vspace{2mm}
 \noindent {\bf Results.} 
In our experiments, we observe that the empirical  performance of affine policy is near-optimal. In particular, the performance is significantly better than the theoretical performance bounds implied in Theorem~\ref{thm:bounded} and Theorem~\ref{thm:gaus}. For instance, Theorem~\ref{thm:bounded} implies that affine policy is a 2-approximation with high probability for i.i.d. random instances from a uniform distribution (see Corollary~\ref{clr:bounded}). However, in our experiments, we observe that the optimality gap for affine policies is at most $4\%$ (i.e. approximation ratio of at most $1.04$). The same observation holds for Gaussian distributions as well Theorem \ref{thm:gaus} gives an approximation bound of $O(\sqrt{\log(mn)})$. We would like to remark that we are not able to report the ratio $r$ for large values of $m$ because the adjustable problem is computationally very challenging and for $m \geq 40$, MIP does not solve within a time limit of $3$ hours for most instances . On the other hand, affine policy scales very well and the average running time is few seconds even for large values of $m$. This demonstrates the power of affine policies that can be computed efficiently and give good approximations for a large class of instances.

 %On the other hand, computing the optimal affine policy over ${\cal U}$ becomes computationally challenging as $m$ increases. For instance, the average running time for computing an optimal affine policy for $m=100$ is around $16$ minutes for the hypersphere uncertainty set, around $36$ minutes for the $3$-norm ball and around $26$ minutes for the $3/2$-norm ball. Whereas, our policy can be computed in less than $1$ second for all instances on average. Figure \ref{fig:time} gives the running time of affine policy for different values of $m$ and as we can observe this running time becomes very high for large values of $m$.

%\iffalse%%%%%%%%%%%%%%%%%%%%%%%%%%%%%%%%%%%%%%%%%%%%%%%%%%%%%%%%%%%%%%%%%%%%%%%%%%%%

\begin{table}[htp]
\begin{subtable}{.5\linewidth}\centering
{\begin{tabular}[t]{|l|l|l|c|c|}
\hline
$m$ & $r_{\sf avg}$ &  $r_{\sf max}$  & $T_{\sf AR}(s)$ & $T_{\sf Aff}(s)$ \\ \hline
10 &1.01 &1.03 & 10.55 &  0.01 \\  \hline
20 & 1.02& 1.04 & 110.57 & 0.23\\  \hline
30 &     1.01          &         1.02         &                 761.21        & 1.29 \\  \hline
50 &** & **& **  & 14.92   \\  \hline
\end{tabular}}
\caption{Uniform}\label{tab:uniform}
\end{subtable}%
\begin{subtable}{.5\linewidth}\centering
{\begin{tabular}[t]{|l|l|l|c|c|}
\hline
$m$ & $r_{\sf avg}$ &  $r_{\sf max}$  & $T_{\sf AR}(s)$ & $T_{\sf Aff}(s)$\\ \hline
10 & 1.00 & 1.03& 12.95& 0.01 \\  \hline
20 &1.01 & 1.03& 217.08& 0.39\\  \hline
30 & 1.01&1.03 &594.15 & 1.15 \\  \hline
50 &** & **&** &  13.87 \\  \hline

\end{tabular}}
\caption{Folded Normal}\label{tab:gaussian}
\end{subtable}
\caption{Comparison on the performance and computation time of affine policy and optimal adjustable policy for uniform and folded normal distributions. For 20 instances, we compute $ z_{\sf{Aff}} (\tilde{\mb B})    /z_{\sf{AR}}  (\tilde{\mb B})         $ and present the average and $\max$ ratios. Here, $T_{\sf AR}(s)$ denotes the running time for the adjustable policy and $T_{\sf Aff}(s)$ denotes the running time for affine policy in seconds. ** Denotes the cases when we set a time limit of 3 hours. These results are obtained using Gurobi 7.0.2 on a 16-core server with 2.93GHz processor and 56GB RAM.}
\label{tab}
\end{table}

\bibliographystyle{informs2014} 
\bibliography{robust} 

%{\singlespace
%{\bibliographystyle{abbrv}
%\bibliography{robust}}}

%\appendix
%\begin{appendix}
\begin{APPENDICES}

\section{Proofs of preliminaries}
\subsection*{ Proof of Lemma \ref{lem:reform}}\label{apx-proofs:lem:reform}
\proof{Proof.}
We have
\begin{align}
z_{\sf AR}({\mb B}) &= \min_{\mb x \geq \mb 0 } \;  \mb{c}^T \mb{x} +  \max_{\mb{h}\in {\cal U}} \; \;  \nonumber
\underset {\mb y \geq \mb 0 }{\underset{  \mb{B}\mb{y}  \; \geq \; \mb{h}-\mb{A}\mb{x}  }\min } \;  \mb d^T \mb y \\ \label{eq:proof:sep}
&= \min_{\mb x \geq \mb 0 } \;  \mb{c}^T \mb{x} +  \max_{\mb{h}\in {\cal U}} \; \; 
\underset {\mb w \geq \mb 0 }{\underset{  \mb B^T \mb w \leq \mb d  }\max } \; (\mb h - \mb A \mb x)^T \mb w \\ \nonumber
&= \min_{\mb x \geq \mb 0 } \;  \mb{c}^T \mb{x} +  \max_{\mb{w}\in {\cal W}} \; - (\mb A \mb x)^T \mb w + \;
\underset {\mb h \geq \mb 0 }{\underset{  \mb R  \mb h \leq \mb r  }\max } \;  \mb h ^T \mb w \\ \nonumber
&= \min_{\mb x \geq \mb 0 } \;  \mb{c}^T \mb{x} +  \max_{\mb{w}\in {\cal W}} \; - (\mb A \mb x)^T \mb w + \;
\underset {\mb \lambda \geq \mb 0 }{\underset{  \mb R ^T \mb \lambda \geq \mb w  }\min } \;  \mb r ^T \mb{\lambda}  \\ \nonumber
&=  z_{\sf d-AR}({\mb B}).
\end{align}
where the second equality holds by taking the dual of the inner minimization problem, the third equality follows from switching the two max, and the fourth one by taking the dual of the second maximization problem.
\hfill
\Halmos
\endproof

\subsection*{Proof of Lemma \ref{lem:berti-goyal}}\label{apx-proofs:lem:berti-goyal}
\proof{Proof.}
We restate the same proof in \cite{BG10} in our setting. First, since the adjustable problem is a relaxation of the affine problem then $z_{\sf d-AR}({\mb B}) \leq  z_{\sf d-Aff}({\mb B}) $.

Now let us prove the other inequality. Consider  ${\cal W}=\{\mb{w}\in{\mathbb R}^m_+\;|\;    \mb B^T \mb w \leq \mb d \}$ which is a simplex. Note that $ \mb 0$ is always an extreme point of the simplex ${\cal W}$ and denote $\mb w^1, \mb w^2, \ldots , \mb w^m$ the remaining $m$ points. In particular, we have for any $ \mb w \in {\cal W}$
$$ \mb w =\sum_{j=1}^m  \alpha_j \mb w^{j} = \mb Q \mb \alpha$$
where $ \sum_{j=1}^m \alpha_j \leq 1$ and $ \mb Q= \left[ \mb w^1 \vert \mb w^2 \vert \ldots \vert \mb w^m  \right]  $. Note that $\mb Q$ is invertible since $\mb w^1, \mb w^2, \ldots , \mb w^m$ are linearly independent. Hence, $ \mb \alpha = \mb Q^{-1} \mb w$.
 Denote $\mb x^*,  \mb \lambda^* ( \mb w )$,$ \mb w \in {\cal W}$, an optimal solution of the adjustable problem \eqref{eq:dual}.  We define the following affine solution $\mb x = \mb x^*$ and for $ \mb w \in {\cal W}$, $$ \mb \lambda(\mb w) = \mb P \mb Q^{-1} \mb w$$ where $$ \mb P= \left[ \mb \lambda^*(\mb w^1) \vert \mb \lambda^*(\mb w^2) \vert \ldots \vert \mb \lambda^*(\mb w^m ) \right]  .$$ In particular, we have 
 $$ \mb \lambda ( \mb w) =  \sum_{j=1}^m \alpha_j \mb \lambda^*(\mb w^j).$$ Let us first check the feasbility of the solution. We have,
 $$ \mb R^T \mb \lambda ( \mb w)  = \sum_{j=1}^m \alpha_j \mb R^T \mb \lambda^*(\mb w^j) \geq \sum_{j=1}^m \alpha_j  \mb w^j = \mb w $$
where the inequality follows from the feasibility of the adjustable solution. Therefore,
\begin{align*}
z_{\sf d-Aff}({\mb B}) &\leq \mb c^T \mb x + \max_{ \mb w \in {\cal W}} \; \; ( -\mb A \mb x )^T \mb w + \mb r^T \mb \lambda ( \mb w)\\
& = \mb c^T \mb x^* + \max_{ \mb \alpha} \; \; ( -\mb A \mb x^* )^T \mb w + \sum_{j=1}^m \alpha_j \mb r^T \mb \lambda^* ( \mb w^j)\\
& = \mb c^T \mb x^* + \max_{ \mb \alpha} \; \;   \sum_{j=1}^m  \alpha_j \left( ( -\mb A \mb x^* )^T \mb w^j + \mb r^T \mb \lambda^* ( \mb w^j) \right)  \\
& \leq \mb c^T \mb x^* + \max_{ \mb w \in {\cal W}}    \left( ( -\mb A \mb x^* )^T \mb w + \mb r^T \mb \lambda^* ( \mb w) \right) \max_{ \mb \alpha} \;  \sum_{j=1}^m \alpha_j   \leq  z_{\sf d-AR}({\mb B}) \\
\end{align*}
where the last inequality holds because $ \sum_{j=1}^m \alpha_j \leq 1$. We conclude that $z_{\sf d-Aff}({\mb B})= z_{\sf d-AR}({\mb B})$.
\hfill
\Halmos
\endproof

\subsection*{Proof of Lemma \ref{lem:inclusion}}\label{apx-proofs:lem:inclusion}
\proof{Proof.}
First the inequality $ z_{\sf d-AR}(\mb B) \leq  z_{\sf d-Aff}(\mb B) $ is straightforward since the adjustable problem\eqref{eq:ar} is a relaxation of the affine problem \eqref{eq:aff}.
On the other hand, since ${\cal W} \subseteq \kappa \cdot {\cal S}$ then,
$$ z_{\sf d-Aff}(\mb B) \leq \kappa  \cdot z_{\sf d-Aff}(\mb B, {\cal S})  $$
where we denote  $z_{\sf d-Aff}(\mb B, {\cal S})$ the dualized affine problem over $ {\cal S}$ (it's the same problem as $z_{\sf d-Aff}(\mb B)$ where we only replace ${\cal W}$ by ${\cal S}$).  Since ${\cal S}$ is a simplex, from Lemma \ref{lem:berti-goyal}, we have $z_{\sf d-Aff}(\mb B, {\cal S})  =z_{\sf d-AR}(\mb B, {\cal S})  $. Moreover, $z_{\sf d-AR}(\mb B, {\cal S})  \leq z_{\sf d-AR}(\mb B) $ because ${\cal S} \subseteq   {\cal W}$. We conclude that 
$$ z_{\sf d-AR}(\mb B) \leq  z_{\sf d-Aff}(\mb B)  \leq \kappa \cdot z_{\sf d-AR}(\mb B).$$
Furthermore, since  $ z_{\sf d-AR}(\mb B)= z_{\sf AR}(\mb B)$ from Lemma \ref{lem:reform} and $ z_{\sf d-Aff}(\mb B)= z_{\sf Aff}(\mb B)$ from Lemma \ref{lem:berti-de Ruiter}, then
$$ z_{\sf AR}(\mb B) \leq  z_{\sf Aff}(\mb B)  \leq \kappa \cdot z_{\sf AR}(\mb B).$$
\hfill
\Halmos
\endproof

\section{Hoeffding's inequality}\label{apx-proofs:Hof-ineq}
 {\bf Hoeffding's inequality}\cite{hoeffding1963probability}.  Let $Z_1, \ldots, Z_n$ be independent bounded random variables with $Z_i \in [a,b]$ for all $i \in [n]$ and denote $Z= \frac{1}{n} \sum_{i=1}^n Z_i$. Therefore, $$ \mathbb{P} ( Z - \mathbb{E}(Z) \leq - \tau ) \leq \exp  \left( \frac{-2n{\tau}^2}{(b-a)^2} \right).$$

\section{LP and MIP formulations for the empirical section} \label{apx-LP-MIP}

{\bf LP formulation for the affine problem.}
The affine problem \eqref{eq:aff} can be formulated as the following LP
\begin{equation}\label{eq:aff:LP}
\begin{aligned}
z_{\sf Aff}({\mb B}) =  \min \; &  \mb{c}^T \mb{x} + z \\
& z- \mb d^T \mb q \geq \mb r^T \mb v \\
& \mb R^T \mb v \geq \mb P^T \mb d \\
& \mb A \mb x + \mb B \mb q \geq \mb V^T \mb r\\
& \mb R^T \mb V \geq \mb I_m - \mb B \mb P \\
& \mb q \geq \mb U^T \mb r \\
& \mb U^T \mb R + \mb P \geq \mb 0 \\
& \mb{x}   \in  {\mathbb R}^{n}_+, \; \mb{v}   \in  {\mathbb R}^{L}_+, \; \mb{U}   \in  {\mathbb R}^{L \times n}_+, \; \mb{V}   \in  {\mathbb R}^{L \times m}_+.  \\
\end{aligned}
\end{equation}

\proof{Proof.}

The affine problem \eqref{eq:aff} can be reformulated as follows
\begin{equation*}
\begin{aligned}
z_{\sf Aff}({\mb B}) = \min_{\mb x} \; & \mb{c}^T \mb{x} + z\\
& z \geq   \mb{d}^T\left( \mb{P}\mb{h}+\mb{q}\right)  \; \; \; \forall \mb h \in {\cal U} \\
& \mb{A}\mb{x} + \mb{B}\left( \mb{P}\mb{h}+\mb{q}\right) \; \geq \; \mb{h}   \; \; \; \forall \mb h \in {\cal U} \\
& \mb{P}\mb{h}+\mb{q} \; \geq \; \mb{0}   \; \; \; \forall \mb h \in {\cal U} \\
& \mb{x}   \in  {\mathbb R}^{n}_+.\\
\end{aligned}
\end{equation*}
We use standard duality techniques to derive formulation \eqref{eq:aff:LP}. The first constraint is equivalent to 
$$ z - \mb d^T \mb q \geq \underset {\mb h \geq \mb 0 }{\underset{  \mb R  \mb h \leq \mb r  }\max } \; \mb{d}^T\mb{P}\mb{h} .$$
By taking the dual of the maximization problem, the constraint is equivalent 
$$ z - \mb d^T \mb q \geq \underset {\mb v \geq \mb 0 }{\underset{   \mb R^T \mb v \geq \mb P^T \mb d  }\min } \;      \mb{r}^T\mb v. $$
We can then drop the min and introduce $\mb v $ as a variable, hence we obtain the following linear constraints
$$ z- \mb d^T \mb q \geq \mb r^T \mb v, \qquad \qquad \mb R^T \mb v \geq \mb P^T \mb d, \qquad \qquad \mb{v}   \in  {\mathbb R}^{L}_+.$$
We use the same technique for the second sets of constraints, i.e.,
$$ \mb{A}\mb{x} + \mb{B}   \mb q   \ \; \geq \;    \underset {\mb h \geq \mb 0 }{\underset{  \mb R  \mb h \leq \mb r  }\max } \;         \mb{h} (  \mb I_m - \mb B \mb P ).  $$
By taking the dual of the maximization problem for each row and dropping the $\min$ we get the following compact formulation of these constraints
$$ \mb A \mb x + \mb B \mb q \geq \mb V^T \mb r, \qquad \qquad \mb R^T \mb V \geq \mb I_m - \mb B \mb P, \qquad \qquad  \mb{V}   \in  {\mathbb R}^{L \times m}_+. $$
Similarly, the last constraint 
$$ \mb{q} \; \geq \;\underset {\mb h \geq \mb 0 }{\underset{  \mb R  \mb h \leq \mb r  }\max } \;      -\mb{P}\mb{h} ,  $$
is equivalent to 
$$ \mb q \geq \mb U^T \mb r, \qquad \mb U^T \mb R + \mb P \geq \mb 0, \qquad \mb{U}   \in  {\mathbb R}^{L \times n}_+ . $$
\hfill
\Halmos
\endproof

\noindent
{\bf MIP formulation for the separation adjustable problem.} The separation problem \eqref{eq:sep} can be formulated as the following MIP
\begin{equation}\label{eq:MIP}
\begin{aligned}
 \max  \; \;  &\sum_{i=1}^m \sum_{j= -\Delta_{\cal W}}^s  \sum_{k= -\Delta_{\cal U}}^s  \frac{1}{2^{j+k}} \cdot \gamma_{ijk}   - (\mb A  \hat{\mb x})^T \mb w  \\
& \mb w = \sum_{i=1}^m \sum_{j= -\Delta_{\cal W}}^s \frac{\beta_{ij}}{2^j} \cdot  \mb e_i \\
& \mb h = \sum_{i=1}^m \sum_{k= -\Delta_{\cal U}}^s \frac{\alpha_{ik}}{2^k} \cdot  \mb e_i \\
& \gamma_{ijk} \leq \beta_{ij} &\forall i \in [m], j \in [ -\Delta_{\cal U} ,s] , k \in [-\Delta_{\cal W},s] \\
& \gamma_{ijk} \leq \alpha_{ik}   &\forall i \in [m], j \in [ -\Delta_{\cal U} ,s] , k \in [-\Delta_{\cal W},s] \\
& \gamma_{ijk} +1 \geq \alpha_{ik} + \beta_{ij}  &\forall i \in [m], j \in [ -\Delta_{\cal U} ,s] , k \in [-\Delta_{\cal W},s] \\
&  \alpha_{ik}, \beta_{ik}, \gamma_{ijk} \in \{0,1\}                        &\forall i \in [m], j \in [ -\Delta_{\cal U} ,s] , k \in [-\Delta_{\cal W},s] \\
& \mb R \mb h \leq \mb r \\
& \mb B^T \mb w \leq \mb d \\
\end{aligned}
\end{equation}
where $s= \lceil{ \log_2 \left( \frac{m}{\epsilon } \right) \rceil}$, $\Delta_{\cal W}$ is an upper bound on any component of $ w \in {\cal W}$, $\Delta_{\cal U}$ is an upper bound on any component of $ h \in {\cal U}$ and $\epsilon$ is the accuracy of the problem.

\proof{Proof.}
The separation problem \eqref{eq:sep} is equivalent to solving the following problem for given $\hat{\mb x}$

$$\underset {\mb w \in \cal W} {\underset{  \mb h \in {\cal U}   }\max } \; {\mb h}^T \mb w  - (\mb A  \hat{\mb x})^T \mb w  $$
The constraints of the above problem are linear and the second term in the objective function is linear as well. So we will focus only on the first term ${\mb h}^T \mb w$  which is a bilinear function and write it in terms of linear constraints and binary variables. Let us write $ \mb h = \sum_{i=1}^m h_i \mb e_i$. For all $i \in [m]$ we digitize the component $\mb h _i$ as follows 
$$h _i= \sum_{k= -\Delta_{\cal U}}^s \frac{\alpha_{ik}}{2^k}$$
where $s= \lceil{ \log_2 \left( \frac{m}{\epsilon } \right) \rceil}$,  $\Delta_{\cal U}$ is an upper bound on any  $ h_i$ and $\alpha_{ik}$ are binary variables. This digitization gives an approximation to $h_i$ within $\frac{\epsilon }{m}$ which translates to an accuracy of $\epsilon$ in the objective function. We have
$$ \mb h = \sum_{i=1}^m \sum_{k= -\Delta_{\cal U}}^s \frac{\alpha_{ik}}{2^k} \cdot  \mb e_i .$$
Similarly, we have
$$ \mb w = \sum_{i=1}^m \sum_{j= -\Delta_{\cal W}}^s \frac{\beta_{ij}}{2^j} \cdot  \mb e_i $$
where $\Delta_{\cal W}$ is an upper bound on any component of $ w \in {\cal W}$. Therefore, the first term in the objective function becomes
$$ \sum_{i=1}^m \sum_{j= -\Delta_{\cal W}}^s  \sum_{k= -\Delta_{\cal U}}^s  \frac{1}{2^{j+k}} \cdot \alpha_{ik}  \beta_{ij} .$$
The final step is to linearize the term $\alpha_{ik}  \beta_{ij}$. We set, $\alpha_{ik}  \beta_{ij}=\gamma_{ijk}$ where again $\gamma_{ijk}$ is a binary variable. Since all the variables here are binary we can express  $\gamma_{ijk}$ using only linear constraints as follows 
$$ \gamma_{ijk} \leq \beta_{ij}$$
$$ \gamma_{ijk} \leq \alpha_{ik}   $$
$$ \gamma_{ijk} +1 \geq \alpha_{ik} + \beta_{ij}  $$
which leads to formulation \eqref{eq:MIP}.

\hfill
\Halmos
\endproof
\noindent

%\end{appendix}
\end{APPENDICES}

\end{document}